\documentclass[12pt]{iopart}
\usepackage{graphicx}

\newcommand{\deriv}[2]{ \frac{\partial#1}{\partial#2}}

\begin{document}

\title{   Diffusion  on a Hypersphere:   Application to the Wright-Fisher  model}

\author{Kishiko Maruyama\\Yoshiaki Itoh}
\address{The Institute of Statistical Mathematics and The Graduate University 
for Advanced Studies, 10-3 Midori-cho, Tachikawa,  Tokyo 190-8562, Japan}

\begin{abstract}
The eigenfunction expansion by Gegenbauer polynomials for the diffusion 
on a hypersphere is transformed into the diffusion for the Wright-Fisher model 
with a particular mutation rate.  We use the Ito calculus considering stochastic 
differential equations. The expansion gives a simple interpretation of the Griffiths
 eigenfunction expansion for the Wright-Fisher model. 
Our representation is useful to simulate the Wright-Fisher model as well as 
Brownian motion on a hypersphere.
\end{abstract}

\pacs{02.50.Ey, 87.10.Ca}

\maketitle

\section{Introduction}
Here we show that the  diffusion on a hypersphere \cite{ca} is transformed 
into the diffusion for the
Wright-Fisher  model with a particular mutation rate \cite{gr1,gr2,gr3}, 
 by  using   the relation, $x_i=y_i^2$, 
where $x_i$'s denote the relative abundance  of alleles 
and $y_i$'s denote the position of a particle  of the diffusion 
on a hypersphere. 
Diffusion on a sphere has been   applied to various problems in   physics, chemistry, 
 mathematics, etc.
   \cite{de,per,yo,br}. 
The Wright-Fisher model is widely used   to study   the random sampling effect  
in population genetics \cite{ew,fi,ki,wr1}.

 We consider   stochastic 
differential equations  in the Ito sense \cite{mc}.   We
 represent 
  the  Wright-Fisher  model for 
$k$ alleles  as well as the   diffusion on   $(k-1)$-dimensional sphere \cite{ma},  
 by using 
  ${k\choose 2}$ 
independent  standard Brownian motion  $b_{ij}(t)$ $i>j$  with  
the skew symmetry  $b_{ij}(t)+b_{ji}(t)=0$ 
\cite{i1,i3, ims, mi}. 
 Our  representation is   useful  to simulate the 
Wright-Fisher  model \cite{mn,ta}.  There are effective simulation methods to simulate 
a 
Brownian motion on a hypersphere for example  \cite{cee}.  Our 
representation is simple and  may be useful to simulate the Brownian 
 motion on a hypersphere. 

\section{ Stochastic differential equation for the  Wright-Fisher  model }

In the Wright-Fisher  model each of the genes of the next generation 
is obtained by a random choice among the genes of the previous 
generation 
and that the whole population changes all at once. 
In Moran's model \cite{mo} it is supposed that there are $2N$ individuals 
each formed from $k$ alleles $A_{1}, A_{2}, ...,A_{k}$, 
and that each instant at which the state of the model may change, 
one individual of the alleles, chosen at random, 
dies and is replaced by a new individual 
which is $A_{i}$ with probability $\displaystyle\frac {k_{i}}{2N}$, 
where $k_{i}$ is the abundance of the allele $A_{i}$. 
It is supposed that the probability of any individual 
"dying" during an interval $(t,t+dt)$ 
and then being replaced by a new individual is $\lambda dt$. 
Hence the mean number of such events in unit time 
is $\lambda$ and the mean length of a generation is $\lambda^{-1}$. 
The following model of interacting particles  \cite{i2} 
is equivalent with the Moran model.

Consider a population of $N$ particles 
each of which is one of $k$ types, $A_{1}, A_{2},..., A_{k}$. 
The types may represent species, alleles, genotypes or other classification.  
We then consider random interactions between particles, 
which are assumed to occur at the rate $\lambda dt$ 
per time interval $(t,t+dt)$ for each particle. 
If a pair of particles of different types $i$ and $j$ interact, 
then after the interaction the both particles are the type $i$ 
with probability 1/2 and the type $j$ with probability 1/2. 
If the type of the interacting  particles are the same, no change occurs.

We can approximate the behavior of our  interacting particle system
by a stochastic differential equation (\ref{wright}). 
In it, the relative abundance of type $i$ increases 
by $c\sqrt {x_{i}(t)x_{j}(t)} \ db_{ij}(t)$ and 
decreases by $c\sqrt {x_{i}(t)x_{j}(t)} \ db_{ij}(t)$ 
with the interaction of the particles of type $j$, 
where $c=\sqrt {\lambda/2N} $. 
Hence our interacting particle system  automatically makes 
the following equation (\ref{eq:collision}), 
which has the genetic drift matrix 
with the elements $c^{2} x_{i}(t)(\delta _{ij}-x_{j}(t))dt$ 
for $i,j=1,2,...,k$ as covariances.

For $i,j=1,2,...,k$, consider
\begin{equation}\label{wright}
  dx_{i}(t)=\sum _{j=1,i\neq j}^{k} 
	c\sqrt {x_{i}(t)x_{j}(t)} db_{ij}(t)	\label{eq:collision}
\end{equation}
with $b_{ij}(t)+b_{ji}(t)=0$,
where $b_{ij}(t)$ $(i>j)$ are mutually independent one-dimensional 
Brownian motion with the mean 0 and the variance $t$ 
\cite{i1,i3,ims,mi}.  
This representation of the Wright-Fisher  model \cite{i1} is applied 
to make computer simulations on models in population genetics \cite{mn,ta}. 
Here we make use of this equation (\ref{eq:collision}) 
to discuss the exact solutions of the diffusion equations 
considering the  diffusion on a unit sphere.

\section{Diffusion  on hypersphere}
 The isotropic diffusion on the $(k-1)$-dimensional unit sphere is given by the
following stochastic differential equation in the  Ito sense. 
\begin{equation}\label{sdesphere}
  dy_{i}(t)=\frac {-c^{2}}{8}(k-1)y_{i}(t)dt+ 
	\frac {c}{2} \sum _{j=1}^{k} y_{j}(t) db_{ij}(t)
\end{equation}
with $b_{ij}(t)+b_{ji}(t)=0$, for $i=1,2,...,k$,
where $b_{ij}(t)$ $(i>j)$ are mutually independent one-dimensional   
Brownian motion with the mean 0 and the variance $t$. 
Let $\mathbf{dy}=(dy_1,...,dy_k)$.
By using the Ito calculus \cite{mc},  we have for the stochastic differential 
\begin{eqnarray}
d (\sum_i y_i(t)^2)= 0.
\end{eqnarray}
Hence starting from a point on the unit hypersphere, the trajectory is on the unit 
hypersphere $\sum_i y_i(t)^2=1$. 
Consider the dot product for a unit vector $\mathbf{l}=(l_1,l_2,...,l_k)$, 
\begin{eqnarray}
\mathbf{l} \cdot \mathbf{dy}= \sum_{i}l_i [\frac{c}{2}\sum_{j} ( y_j (t) db_{ij}(t))
+\frac {-c^{2}}{8}(k-1)y_{i}(t)dt].
\end{eqnarray}
On  the tangent hyperplane at $\mathbf{y}=(y_1,y_2,...,y_k)$ of the hypersphere, 
 for   $\mathbf{l}$ which is perpendicular to $\mathbf{y}$,  we have  
 \begin{eqnarray}
\mathbf{l} \cdot \mathbf{dy} =\frac{c}{2} \sum_{i}l_i (\sum_{j}  y_j (t) db_{ij}(t)).
 \end{eqnarray}
 We  have  the expectation 
 \begin{eqnarray}
E(\mathbf{l} \cdot \mathbf{dy})=0.
\end{eqnarray}
We have 
\begin{eqnarray}\label{var}
Var(\mathbf{l} \cdot \mathbf{dy})=(\frac{c}{2})^2\sum_{i>j}(l_iy_j-l_jy_i)^2dt=(\frac{c}{2})^2dt,
\end{eqnarray}
since
\begin{eqnarray}
\nonumber&&\sum_{i>j}(l_iy_j-y_il_j)^2=\sum_{i> j} [(l_iy_j)^2+(y_il_j)^2-2l_iy_jy_il_j]+(\sum_i \l_iy_i)^2\\
&& 
=\sum_{i>j} [(l_iy_j)^2+(y_il_j)^2]+\sum_i (l_iy_i)^2=(\sum_{i}l_i^2)(\sum_{i}y_i^2)=1.
\end{eqnarray}
$Var(\mathbf{l} \cdot \mathbf{dy})$  does not depend on the direction  $\mathbf{l}$ 
at the tangential 
point $\mathbf{y}$ of  the tangential plane.  Hence we see the diffusion is isotropic 
and we have the Fokker-Planck equation on the hypersphere for the 
Laplace-Beltrami operator  \cite{ca} in space $S_{k-1}$,  
\begin{eqnarray}\label{rh}
\frac{\partial}{ \partial t} \rho(\mathbf{y},t|\mathbf{y'},0) =D\triangle^{S_{k-1}} 
 \rho(\mathbf{y},t|\mathbf{y'},0) 
\end{eqnarray}
where $D=\displaystyle\frac{c^{2}}{8}$ from (\ref{var}).
The operator  $\triangle^{S_{k-1}} $  is for $k - 1$ spherical 
 coordinates $\theta _1 , \theta _2 , \dots , \theta _{k-1}$ \, 
of the unit
vector $y = (y_1, . . . , y_k)$ of $S_{k-1}$  defined by the relations
\begin{eqnarray}
y_1 &&= \sin\theta_{k-1} . . . \sin \theta_2 \sin  \theta_1\\
y_2 &&= \sin  \theta_{k-1} . . . \sin  \theta_2 \cos  \theta_1\\
&&\vdots\\
 y_{k-1} &&= \sin  \theta_{k-1} \cos  \theta_{k-2} \\
 y_k &&= \cos  \theta_{k-1}, 
 \end{eqnarray}
where $0 \leq  \theta_1 < 2 \pi$ and $0 \leq \theta_i < \pi$ for $i \ne 1$.  
 The integration measure in 
$S_{k-1}$ is  defined as
$dy=\frac{1}{A_{k-1}}
\sin^{k-1} \theta_{k-1} . . . \sin  \theta_1 d \theta_1 . . . d \theta_{k-1}$ , 
where $A_{k-1} = 2\pi^{k/2}/\Gamma(k/2)$ is the surface of the sphere $S_{k-1}$.  
With this normalization 
$\int_{S_{k-1}}dy= 1$.

The solution \cite{ca},  the transition probability 
density for the isotropic diffusion on hypersphere, is  obtained   by using 
the definitions and notations of the chapter IX of the book \cite{vi}. 
The Gegenbauer polynomial $C^p_L$ is defined as the coefficient of $h^L$ 
in the power-series expansion of the
function
\begin{eqnarray}
(1 - 2th + h^2)^{-p} =
\sum^{\infty}_{L=
0}
C^p_L(t)h^L . 
\end{eqnarray}
The transition probability density  from $\mathbf{y'}$   (with  $n - 1$ spherical 
 coordinates $\theta' _1, \dots , \theta' _{n-1}$)  at $0$ 
to $\mathbf{y}$ (with  $\theta _1, \dots , \theta _{n-1}$ )
at $t>0$    is given \cite{ca} by 
\begin{eqnarray}\label{rho}
&&\rho(\mathbf{y},t|\mathbf{y'},0) \nonumber\\&&=
\frac{1}{
A_{k-1}}
\sum^{\infty}_{L=0}\frac{2L + k - 2}{k - 2}
C^{k/2-1}_L (\mathbf{y}\cdot \mathbf{y'})
  \exp(
-DL(L + k - 2) t), 
\end{eqnarray}
for the dot product $\mathbf{y}\cdot \mathbf{y'}=\sum_i y_iy'_i$,  where 
\begin{eqnarray}\label{gegen}
C_L^{(k/2-1)}(z)=\sum_{j=0}^{\lfloor L/2 \rfloor}(-1)^j
\frac{\Gamma(L-j+k/2-1)}
{\Gamma(k/2-1)j!(L-2j)!}(2z)^{L-2j}.
\end{eqnarray}

\section{Solution for the 
Wright-Fisher 
model with parent-independent mutation}

Let $A_{1},...,A_{k}$ denote $k$ allele types in a population. 
The general neutral alleles model has a probability $u_{ij}$ 
of a mutation from $A_{i}$ to $A_{j}$. 
 An exact solution is known for  the case of the 0 mutation rate  case 
in three dimensions \cite{ki}.
Assuming parent-independent mutation \cite{gr1}, $u_{ij}=u_{j}$.  
an  expansion of the transition probability density \cite{gr1},   
from  the initial relative abudances  
$\mathbf{x'}=(x'_{1},..., x'_{k})$ to  the  relative abundances 
$\mathbf{x}=(x_{1}, ..., x_{k})$, 
is obtained for  the following Fokker-Planck equation 
(\ref{eq:fp3}), 
\begin{equation}
  \deriv {p(\mathbf{x},t\mid \mathbf{x'},0)}{t}=\sum_{i=1}^{k-1}\deriv { }{x_{i}}(\frac {1}{2}
 	\sum_{i=1}^{k-1}\deriv {x_{i}(\delta_{ij}-x_{j})}{x_{j}}-\frac {1}{2}M_{i})p(\mathbf{x},t\mid \mathbf{x'},0)
							\label{eq:fp3}
\end{equation}
where $M_{i}=\varepsilon_{i}-\mu x_{i}$, 
$\mu=\sum_{i=1}^{k} \varepsilon_{i}$. 
For $\varepsilon_{i}>0, i=1,2,...,k$, 
the stationary density  \cite{wr2, wa} is given by  
\begin{equation} 
  \Gamma(\mu)\prod _{i=1}^{k}x_{i}^{\varepsilon_i-1}/
  	\Gamma(\varepsilon_i).			
  	\label{eq:stationary}
\end{equation}
The solution of equation  (\ref{eq:fp3}) is given by an expansion in orthogonal
polynomials  \cite{gr1}.   For the 
case $\varepsilon_i=\varepsilon$, $i=1,...,k$,    we have the solution \cite{gr1}, 
\begin{eqnarray}\label{Grifsol}
 \nonumber&&p(\mathbf{x},t\mid \mathbf{x'},0)\\ =&&\Gamma(\mu)\prod_{i=1}^{k}\frac{x_{i}^{\varepsilon-1}}{
	\Gamma(\varepsilon)}(\sum_{n=0}^{\infty}\exp\{-\frac{1}{2}n(n-1)t-\frac{1}{2}\mu nt\}
	Q_{n}(\mathbf{x},\mathbf{x'}))
							\label{griffiths}
\end{eqnarray}
for $x_{j}>0, j=1, 2,...,k$, 
where 
 $Q_{n}$ is given  by 
\begin{eqnarray}\label{Qn}
Q_n(\mathbf{x},  \mathbf{x'})=(\mu +2n-1)(n!)^{-1}\sum_{m=0}^n (-1)^{n-m}{n\choose m}
(\mu+m)_{(n-1)}\xi_m
\end{eqnarray}
with
\begin{eqnarray}
\xi_m=\mu_{(m)} \Gamma(\varepsilon)^k\sum_{l_1+...+l_k=m, 0\leq l_j, j=1,...,k}
\frac{m!}{l_1!...,l_k!}\Pi_{j=1}^{k}\frac{(x_j x'_j)^{l_j}}{\Gamma(l_j+\varepsilon)},
\end{eqnarray}
for the notation $\alpha_{(m)}=\alpha( \alpha+1)\cdots( \alpha +m-1)$.  
We give the particular cases as, 
\begin{eqnarray}
Q_0(\mathbf{x},\mathbf{x'})&&=1,\\
Q_1(\mathbf{x},\mathbf{x'})&&=(\mu+1)(\xi_1-1),\\
Q_2(\mathbf{x},\mathbf{x'})&&=\frac{1}{2}(\mu+3)[(\mu+2)\xi_2-2(\mu+1)\xi_1+\mu)],
\end{eqnarray}
with 
\begin{eqnarray}
\xi_1&&=\mu\sum_{j=1}^k\frac{x_j x'_j}{\varepsilon}\\
\xi_2&&=\mu(\mu+1)\sum_{j=1}^k(\frac{(x_j x'_j)^2}{\varepsilon(\varepsilon+1)}
+2\sum_{i<j}\frac{x_i x_j x'_i x'_j}{\varepsilon^2}).
\end{eqnarray}
For $i=1,...,k$,  putting $y_{i}^{2}=x_{i}$,  the isotropic diffusion (\ref{sdesphere}) 
on $(k-1)$-dimensional 
sphere is transformed to
\begin{equation}
  dx_{i}(t)=\frac{c^{2}}{4}(1-kx_{i}(t))dt+
	\sum_{j=1}^{k}c\sqrt{x_{i}(t)x_{j}(t)}db_{ij}(t). 
							\label{eq:isotropic}
\end{equation}
The diffusion process governed 
by equation (\ref{eq:fp3}) is expressed as, 
\begin{equation}
  dx_{i}(t)=M_{i} dt+
	\sum_{j=1}^{k}\sqrt{x_{i}(t)x_{j}(t)}db_{ij}(t). 
							\label{eq:genetic}
\end{equation}
Thus the corresponding Fokker-Planck equation of Eq.(\ref{eq:isotropic}) 
is a special case of the equation (\ref{eq:fp3}), 
putting $\varepsilon_{i}=1/2$ for all $i$ 
and taking $c=1$.

The transition probability density  $p(\mathbf{x},t\mid \mathbf{x'},0)$ 
for  equation  (\ref{eq:fp3}) with $\varepsilon_i=1/2$, $i=,...,k$
  on the relative abundances $\mathbf{x}=(x_1,...,x_k)$   starting from 
$\mathbf{x'}=(x'_1....,x'_k)$  
 is represented by using 
the transition probability density $\rho(  \mathbf{y},t\mid \mathbf{y'} ,0)$,  
for   equation (\ref{rho}),  on 
$\mathbf{y}=(y_1,...,y_k)$ with spherical coordinates  $\theta_1, ...,\theta_k$  of 
$(k-1)$-dimensional sphere starting from $\mathbf{y'}=(y'_1,...,y'_k)$   
with spherical coordinates $\theta'_1, ...,\theta'_k$, 
as 
\begin{eqnarray} 
&&p(\mathbf{x},t\mid \mathbf{x'},0)
\left|\frac{\partial(x_1,...,x_{k-1})}{\partial(\theta_{1},...,\theta_{k-1})}\right|
d\theta_1...d\theta_{k-1}\\&&=
         \sum_{y^{'2}_i=x'_i, 0\leq y'_i,   y^2_i=x_i,  i=1,...,k } 
	\rho(  \mathbf{y},t\mid \mathbf{y'} ,0)
	\left|\frac{\partial(y_1,...,y_{k-1})}{\partial(\theta_{1},...,\theta_{k-1})}\right|
d\theta_1...d\theta_{k-1}.
\label{main}
\end{eqnarray}
We have 
\begin{eqnarray}\label{jacobi3}
\left|\frac{\partial(x_1,...,x_{n-1})}{\partial(\theta_{1},...,\theta_{n-1})}\right|=
2^{k-1}(\prod_{i=1}^{k} y_{i})
\left|\frac{\partial(y_1,...,y_{k-1})}{\partial(\theta_{1},...,\theta_{k-1})}\right|.
\end{eqnarray}
Considering equation  (\ref{rho}), (\ref{surface}), (\ref{gegen}), (\ref{jacobi3}) we have 
\begin{eqnarray}\label{solution}
\nonumber&&  p(\mathbf{x},t\mid \mathbf{x'},0)\\\nonumber
 &=&\frac{\Gamma (k/2)}
 {\pi^{k/2}} (\prod_{i=1}^{k} x_{i}^{-1/2} )~~
2^{-k}\sum_{y^{'2}_i=x'_i, 0\leq y'_i,   y^2_i=x_i,  i=1,...,k } \\&&
\sum^{\infty}_{L=0}\frac{2L + k - 2}{k - 2}
C^{k/2-1}_L (\mathbf{y}\cdot \mathbf{y'})
  \exp(
-DL(L + k - 2) t). 
\end{eqnarray}

Thus we proved the solution  (\ref{solution})  is  equivalent
 with the solution  (\ref{griffiths}) for the case $\varepsilon=1/2$.  
Let us observe the two forms of the solution.
Consider the eigen values of the two diffusion equations (\ref{rh}) and (\ref{eq:fp3}). 
The  Caillol solution \cite{ca}   of (\ref{rh}),  expanded by 
the Gegenbauer   polynomials,   are applied  in equation (\ref{solution}).   
The Gegenbauer   polynomials for odd $L$ 
in  equation  (\ref{solution}) cancel with each other  as we can see from 
equation (\ref{gegen}).  Hence there is no  term for odd $L$ in equation 
 (\ref{solution}), \cite{ma}.
Put $L=2n$, and    take $D=1/8$  ($c=1$ in (\ref{var})),  
for  particular mutation rate $\varepsilon_i=1/2,~~ i=1,...,k$.  
Then  the exponential  functions $\exp(
-DL(L + k - 2) t)$ in equation (\ref{solution}) coincide with those  of  Griffiths 
solution  (\ref {griffiths}).  The above $Q_0(\mathbf{x},\mathbf{x'})$, 
$Q_1(\mathbf{x},\mathbf{x'})$, and $Q_2(\mathbf{x},\mathbf{x'})$ in the 
 solution   (\ref {griffiths})  are  
 obtained from 
the terms in 
the expansion 
(\ref {solution}) by using $x_i=y^2_i$ and $x'_i=y'^2_i$ for $i=1,...,k$, as 
we can calculate  easily.  
Also we see 
\begin{eqnarray}
\frac{\Gamma (k/2)}
 {\pi^{k/2}} (\prod_{i=1}^{k} x_{i}^{-1/2} )
 =\Gamma(\mu)\prod_{i=1}^{k}\frac{x_{i}^{\varepsilon-1}}{
	\Gamma(\varepsilon)}.
	\end{eqnarray}
Our  study  naturally  connects the diffusion on  a hypersphere with 
the Wright-Fisher model for  a particular mutation rate   of the  
  parent-independent mutation.

Remark 1. 
The connection between diffusion on a hypersphere and the 1-dimensional
Wright-Fisher diffusion has been derived before (p338 in \cite{kt} ).  
The equation (13.33) in \cite{kt}  corresponds  to equation (\ref{rho}). 
They derive the spectral expansion for a wider range of parameters. 
To  extend this argument to the multidimensional case is one of our next problems.

Remark 2.  
The eigenfunctions in our equation (\ref{solution}) and the eigen function 
 (\ref{Qn}) in the 
Griffiths solution for the particular mutation rate look different,  although they 
  are identical.  
We do not  have the algebraic proof in our present work and wish to  have it 
  in our next work.   
The form of the polynomials in our  expansion means that they are 
the reproducing kernel polynomials for orthogonal polynomials on the 
Dirichlet measure \cite{gs},   because all pairs of polynomials of the same degree are added
to form the coefficient of a given eigenvalue . 
 
\vspace{0.2cm}
{\it Acknowledgement.} The authors thank  
Hidenori  Tachida for the references  in population genetics.
The authors thank    Shinichi  Kotani, Tohru  Ogawa 
and Geoffrey  A. Watterson for helpful discussion and suggestion
on the symmetry of  special  functions.   The authors thank the reviewers 
for  helpful comments and suggestions.

\end{document}